\RequirePackage[reqno]{amsmath}
\documentclass[preprint]{amsart}

\usepackage{amssymb}
\usepackage{empheq}
\usepackage{comment}
\usepackage{amsthm}
\usepackage{enumerate}
\usepackage[top=1in, left=1in, right=0.9in, bottom=1in]{geometry}

\usepackage{graphicx}
\usepackage{epsfig}
\usepackage{comment}

\newtheorem{theorem}{Theorem}[section]
\newtheorem{lemma}[theorem]{Lemma}
\newtheorem{definition}[theorem]{Definition}
\newtheorem{proposition}[theorem]{Proposition}

\def\R{{\mathbb R}}

\def\S{{\mathbb S}}

\newcommand{\cF}{\mathcal{F}}

\newcommand{\ds}{\displaystyle{}}

\newcommand{\cp}{\operatorname{cap}}
\newcommand{\supp}{\operatorname{supp}}

\begin{document}

\title{Weighted Energy Problem on the Unit Sphere}

\author[Mykhailo Bilogliadov]{Mykhailo Bilogliadov}
\email{mykhail@okstate.edu}
\address{Department of Mathematics, Oklahoma State University, Stillwater, OK 74078, U.S.A.}


\begin{abstract}
We consider the minimal energy problem on the unit sphere $\S^2$ in the Euclidean space $\R^3$ 
immersed in an external field $Q$, where the charges are assumed to interact via Newtonian 
potential $1/r$, $r$ being the Euclidean distance. The problem is solved by finding the support of the 
extremal measure, and obtaining an explicit expression for the equilibrium density. We then apply our
results to an external field generated by a point charge, and to a quadratic external field.


\end{abstract}

\maketitle

\medskip
\medskip
\noindent

\section{Introduction}

Let $\S^2:=\{x\in\R^3: |x|=1\}$ be the unit sphere in $\R^3$, where $|\cdot|$ is the Euclidean distance. Given a compact set $E\subset\S^2$, consider the class $\mathcal{M}(E)$ of unit positive Borel measures supported on $E$. The Newtonian potential and Newtonian energy of a measure $\mu\in\mathcal{M}(E)$ are defined respectively as 
$$U^{\mu}(x):=\int \frac {1} {|x-y|} \, d\mu(y), \quad I(\mu):=\iint \frac {1} {|x-y|} \, d\mu(x)d\mu(y).$$
Let 
$$W(E):=\inf\{I(\mu): \mu\in\mathcal{M}(E)\}.$$
Define the Newtonian capacity of $E$ as $\cp(E):=1/W(E)$. We say that a property holds quasi-everywhere (q.e.), if the exceptional set has Newtonian capacity zero. When $\cp(E)>0$, there is a unique $\mu_{E}$ such that $I(\mu_E)=W(E)$. Such $\mu_E$ is called the Newtonian equilibrium measure for $E$.

An external field is defined as a non-negative lower-semicontinuous function $Q: E \rightarrow [0,\infty]$, such that $Q(x)<\infty$ on a set of positive Lebesgue surface measure. The weighted energy associated with $Q(x)$ is then defined by $$I_Q(\mu):=I(\mu)+2\int Q(x)d\mu(x).$$

The energy problem on $\S^2$ in the presence of the external field $Q(x)$ refers to the minimal quantity 
$$V_Q:=\inf\{I_Q(\mu): \mu\in\mathcal{M}(E)\}.$$ 
A measure $\mu_Q \in\mathcal{M}(E)$ such that $I_Q(\mu_Q)=V_Q$ is called an extremal (or positive Newtonian equilibrium) measure associated with $Q(x)$.

The potential $U^{\mu_Q}$ of the measure $\mu_Q$ satisfies the Gauss variational inequalities
\begin{empheq}{align}\label{var1}
U^{\mu_Q}(x)+Q(x)&\geq F_Q \quad \text{q.e. on}\ \ E,\\\label{var2}
U^{\mu_Q}(x)+Q(x)& \leq F_Q \quad \text{for all}\ \ x\in S_Q,
\end{empheq}
where $F_Q:=V_Q-\ds{ \int Q(x)\,d\mu_Q(x)}$, and $S_Q:=\supp{\mu_Q}$ (see Theorem 10.3 in \cite{mizuta} or Proposition 3 in \cite{bds2}, and also book \cite{bhs}). We remark that for continuous external fields, the inequality in $(\ref{var1})$ holds everywhere, which implies that equality holds in $(\ref{var2})$.

The energy problems with external fields on the sphere were a subject of investigation by 
Brauchart, Dragnev and Saff, see \cite{bds1}--\cite{bds3}, Dragnev and Saff \cite{ddss}. A comprehensive treatment of the subject is contained in a forthcoming book \cite{bhs}.  In paper \cite{bds1}, the authors 
considered the case of the unit sphere $\S^{d-1}\subset \R^d$, $d\geq3$, immersed in an external field, 
generated by a positive point charge, assuming that the point charges interact via the Riesz $s$-potentials $1/r^s$, 
where $r$ is the Euclidean distance between the point charges.  For the case when an external field is 
generated by a positive point charge placed on the positive polar semi-axis outside of the sphere $\S^{d-1}$, and 
more generally, by a line charge, they obtained an equation describing the support of the extremal measure, and also found 
an explicit representation for the equilibrium density. The extensions of the results of \cite{bds1} to the case
of the positive point charge placed inside of the sphere $\S^{d-1}$, as well as to the external fields generated 
by the negative point charges, were considered in \cite{bds2} and \cite{bds3}.

The main aim of this paper is to provide a solution to the weighted energy problem on the sphere 
$\S^2\subset\R^3$, immersed in a general external field, possessing rotational symmetry with respect to the polar 
axis, when support is a spherical cap. It is assumed 
that the charges interact via Newtonian potential $1/r$, with $r$ denoting the Euclidean distance 
between the point charges. We obtain an equation that describes the support of the 
extremal measure, while also giving an explicit expression for the equilibrium density. We then apply our results
to an external field produced by a positive point charge placed on the positive polar semi-axis, thus, in particular, 
recovering the corresponding results of \cite{bds1}, \cite{bds3} for the case of Newtonian potential in $\R^3$. An 
application to an external field, given by a quadratic polynomial, is also considered.

We remark that Theorems \ref{transceqsupp}  and \ref{eqdgct} below recover the corresponding statements first obtained in \cite{bds1}, \cite{bds3}, for the case of Newtonian potential in $\R^3$. However, our results contain no restriction on the hight of the point charge, as long as it is placed on the positive half of the polar axis, and the support is a spherical cap centered at the South Pole. This allows for a unified treatment of the situations when the point charge is placed inside or outside of the sphere $\S^2$, as well as at the North Pole. Note that in \cite{bds1} and \cite{bds3}, an important ingredient of the proofs was the Kelvin transformation with center at the location of the point charge. That dictated an obvious restriction for the point charge to be placed not on the surface of the sphere.

Let $(r,\theta,\varphi)$ be spherical coordinates in $\R^3$, with $0\leq r \leq \infty$, $0\leq \varphi \leq \pi$, 
and $0\leq \theta\leq 2\pi$. A spherical cap on the sphere $\S^2$, centered at the North Pole, is defined 
via an angle $\alpha$, $0 < \alpha \leq \pi$, as
\begin{equation*} 
C_{N,\alpha}:= \{ (r,\theta,\varphi): r=1,\ \ 0\leq \varphi \leq \alpha,\ \  0\leq \theta \leq 2\pi  \}.
\end{equation*}
Similarly, a spherical cap centered at the South Pole, is defined in terms of an angle $\alpha$, $0 < \alpha \leq \pi$, as
\begin{equation*} 
C_{S,\alpha}:= \{ (r,\theta,\varphi): r=1,\ \ \alpha \leq \varphi \leq \pi,\ \  0\leq \theta \leq 2\pi  \}.
\end{equation*}

We begin by recording sufficient conditions on an external field $Q$, that guarantee that the extremal support 
$S_Q$ of the equilibrium measure $\mu_Q$ is a spherical cap $C_{S,\alpha}$, centered at the South Pole.
The following proposition is a consequence of a more general statement, proved in \cite{bds1}.
\begin{proposition}\label{scssp}
Let the external field $Q$ be rotationally invariant about the polar axis $x_3$, that is $Q(x)=Q(x_3)$, where 
$x=(x_1,x_2,x_3)\in\S^2$. Suppose that $Q(x_3)$ is an increasing convex function on $[-1,1]$. Then the support of the 
extremal measure $\mu_Q$ is a spherical cap centered at the South Pole, that is $S_Q=C_{S,\alpha}$.
\end{proposition}

The following result is an important step towards the recovery of the equilibrium measure.
\begin{theorem}\label{theo1}
Suppose that $Q\in C^2(N)$, where $N$ is an open neighborhood of $S_Q=C_{S,\alpha}$ on the sphere $\S^2$. 
Then the equilibrium measure $\mu_Q$ is absolutely continuous with respect to the Lebesgue 
surface measure, with a locally bounded density, that is $d\mu_Q =  f(\varphi)\,\sin\varphi\, d\theta\,d\varphi$, 
where $f\in L^\infty([\alpha,\pi])$.
\end{theorem}

The support $S_Q$ is a main ingredient in determining the equilibrium measure $\mu_Q$ itself. Indeed, if 
$S_Q$ is known, the equilibrium measure $\mu_Q$ can be recovered by solving the singular integral equation
\begin{equation}\label{ieec}
\int \frac{1}{|x-y|}\,d\mu(y)+Q(x)=F_Q, \quad x\in S_Q,
\end{equation}
where $F_Q$ is a constant (see (\ref{var2})).

We solve this equation using the results of \cite{col}, and thus obtain the following two theorems, that describe explicitly the 
equilibrium density when support $S_Q$ is either a spherical cap $C_{N,\alpha}$, centered at the North Pole, or $C_{S,\alpha}$, a spherical
cap centered at the South Pole.
\begin{theorem}\label{theo2}
Let $Q\in C^2(N)$, where $N$ is an open neighborhood of $S_Q$ in $\S^2$. Assume that $S_Q$ is a spherical
cap $C_{N,\alpha}$ centered at the North Pole, with $0 < \alpha \leq \pi$. Let
\begin{equation}\label{auxf}
F(\varphi) = \frac {2} {\pi} \frac{1}{\sin\varphi} \, \frac {d} {d\varphi} \int_{\varphi}^\alpha 
\frac{ g(t) \sin t\, dt} {\sqrt{\cos\varphi - \cos t}},
\end{equation}
where
\begin{equation}\label{auxs}
g(t) = \frac {1} {4\pi}\, \frac {d} {dt} \int_0^t \frac{Q(\xi) \sin \xi  \,d\xi} {\sqrt{\cos \xi - \cos t}}.
\end{equation}
Then the density $f$ of the equilibrium measure $\mu_Q$ is given by
\begin{equation}\label{eqd}
f(\varphi) = \frac {F_Q} {4\pi} \left\{ 1 + \frac {2} {\pi} \left[  \sqrt{ \frac {1+\cos \alpha} {\cos\varphi - \cos \alpha} } -\tan^{-1} \sqrt{ \frac {1+\cos \alpha} {\cos\varphi - \cos \alpha}}   \right]    \right\} + F(\varphi),
\end{equation}
where $0\leq \varphi \leq \alpha $, and the constant $F_Q$ is uniquely defined by
\begin{equation}\label{dcf}
F_Q = \frac {\pi} {\sin \alpha + \alpha} \left\{ 1 - 2\pi \int_0^\alpha F(\varphi) \sin\varphi\,d\varphi \right\}.
\end{equation}
\end{theorem}
\noindent An analogous statement for the support being a spherical cap $C_{S,\alpha}$ centered at the South Pole, 
is of the following nature.
\begin{theorem}\label{theo3}
Let $Q\in C^2(N)$, where $N$ is an open neighborhood of $S_Q$ in $\S^2$. Assume that $S_Q$ is a spherical
cap $C_{S,\alpha}$, centered at the South Pole, with $0 < \alpha \leq \pi$. Let
\begin{equation}\label{auxfs}
F(\varphi) = \frac {2} {\pi} \frac {1} {\sin\varphi} \, \frac {d} {d\varphi} \int_\alpha^{\varphi} 
\frac{ g(t) \sin t \,dt} {\sqrt{\cos t - \cos\varphi}},
\end{equation}
where
\begin{equation}\label{auxss}
g(t) = \frac {1} {4\pi}\, \frac {d} {dt} \int_t^\pi \frac{ Q(\xi) \sin \xi \,d\xi} {\sqrt{\cos t - \cos \xi}}.
\end{equation}
Then the density $f$ of the equilibrium measure $\mu_Q$ is given by
\begin{equation}\label{eqds}
f(\varphi) = \frac {F_Q} {4\pi} \left\{ 1 + \frac {2} {\pi} \left[  \sqrt{ \frac {1- \cos \alpha} {\cos \alpha - \cos\varphi} } -\tan^{-1} \sqrt{\frac {1 - \cos \alpha} {\cos \alpha - \cos\varphi} }   \right]    \right\} + F(\varphi),
\end{equation}
where $\alpha \leq \varphi \leq \pi $, and the constant $F_Q$ is uniquely defined by
\begin{equation}\label{dcfs}
F_Q = \frac {\pi} {\sin \alpha + \pi - \alpha} \left\{ 1 - 2\pi \int_\alpha^\pi F(\varphi) \sin\varphi\,d\varphi \right\}.
\end{equation}
\end{theorem}


\section{Applications}
\noindent We begin by considering the case of no external field, when the support is a spherical 
cap centered at the South Pole, that is $S_Q=C_{S,\alpha}$, $Q=0$. The equilibrium density for 
the spherical cap $C_{N,\alpha}$, centered at the North Pole, for the case of no external field, was 
first found by Collins \cite{col}. Note that the equilibrium measure obtained in \cite{col}, does not
have mass 1. The equilibrium measure for the spherical cap centered at the South Pole $C_{S,\alpha}$, 
was also found in \cite{ddss}, and in \cite{bds1}.

\begin{proposition}
The density of the equilibrium measure of a spherical cap $C_{S,\alpha}$ with no external field is
\begin{equation}\label{emnf}
f(\varphi) = \frac {1} {4 (\pi - \alpha + \sin \alpha)}  \bigg\{ 1 + \frac {2} {\pi} \bigg[  \sqrt{ \frac {1- \cos \alpha} {\cos \alpha - \cos\varphi} }  - \tan^{-1} \sqrt{ \frac {1 - \cos \alpha} {\cos \alpha - \cos\varphi} }   \bigg]    \bigg\}, \quad \alpha \leq \varphi \leq \pi.
\end{equation}
The capacity of $C_{S,\alpha}$ is
\begin{equation}\label{capscs}
\cp(C_{S,\alpha}) = \frac {1} {\pi} (\pi - \alpha + \sin \alpha).
\end{equation}
\end{proposition}

Suppose now that the sphere $\S^2$ is immersed in an external field $Q$, that satisfies the conditions of 
Proposition \ref{scssp}. Then the support $S_Q$ of the weighted equilibrium measure $\mu_Q$ will be 
a spherical cap $C_{S,\alpha}$, centered at the South Pole. One of the ways to determine the angle $\alpha$, 
which defines the extremal support $C_{S,\alpha}$, is via the Newtonian analog of the Mhaskar-Saff $\cF$-functional. Originally 
introduced for the case of logarithmic potential \cite{st}, it was extended by Brauchart, Dragnev and Saff \cite{bds1} 
to the case of the Riesz potentials on the sphere $\mathbb S^{d-1}$, for $d\geq 3$. Here we adapt the approach 
from \cite{bds1}. We begin with a definition.
\begin{definition}\label{msfdef}
The  $\cF$-functional of a compact subset $E\subset \S^2$ of positive (Newtonian) capacity is defined as
\begin{equation}\label{msff}
\cF(E) := W(E) + \int Q(x)\,d\mu_E(x), 
\end{equation}
where $W(E)$ is the Newtonian energy of the compact $E$ and $\mu_E$ is the equilibrium measure (with no external field) on $E$.
\end{definition}
The main objective of introducing the $\cF$-functional is its following extremal property, proved in \cite{bds1} for the general Riesz potentials.
\begin{proposition}\label{msfmp}
Let $Q$ be an external field on $\S^2$. Then $\cF$-functional is minimized for $S_Q=\supp(\mu_Q)$. 
\end{proposition}
\noindent For $E=C_{S,\alpha}$, the functional $\cF$ takes the form
\begin{empheq}{align}\label{mhsf}
\cF(C_{S,\alpha})  = \frac {\pi} {2(\pi-\alpha+\sin \alpha)} \bigg\{2  & + \int_\alpha^\pi Q(\varphi) \sin\varphi\, d\varphi \\\nonumber
							& + \frac {2} {\pi} \int_\alpha^\pi Q(\varphi) \sqrt{ \frac {1-\cos \alpha} {\cos \alpha - \cos\varphi} } \sin\varphi\, d\varphi \\\nonumber
			& - \frac {2} {\pi} \int_\alpha^\pi Q(\varphi) \tan^{-1} \sqrt{ \frac {1-\cos \alpha} {\cos \alpha - \cos\varphi} } \sin\varphi\, d\varphi \bigg\}.
\end{empheq}

\noindent Consider the case of an external field $Q$, generated by a positive point charge of magnitude $q$, hovering over the sphere 
$\S^2$ on the $x_3$-axis at a distance $h$ above the origin,
\begin{equation}\label{ptch}
Q(x) = \frac {q} {|x-hp|}, \quad q>0, \quad p=(0,0,x_3),\ \ x\in\S^2.
\end{equation}
It is clear that $Q(x)=Q(x_3)=q(1+h^2-2hx_3)^{-1/2}$ is a convex increasing function on $[-1,1]$, and is symmetric 
with respect to rotations about the polar axis $x_3$. From Proposition \ref{scssp} it then follows that the support of 
the equilibrium measure $\mu_Q$ will be a spherical cap, centered at the South Pole, that is $S_Q=C_{S,\alpha}$, for some $\alpha\in[0,\pi)$. 

The first step is to compute the $\cF$-functional for such external field $Q$.
\begin{proposition}\label{msfptch}
The $\cF$-functional for the spherical cap $C_{S,\alpha}$ in the case of $Q$ as in $(\ref{ptch})$, has the
form
\begin{equation}\label{msfexplf}
\cF(C_{S,\alpha}) = \frac {\pi} {\pi-\alpha+\sin \alpha} \bigg\{ 1 + \frac {q(h+1)}{2h} \bigg(1 - \frac {\alpha} {\pi} \bigg) - \frac {q(h-1)}{\pi h} \tan^{-1}\bigg( \cot\bigg(\frac{\alpha}{2}\bigg) \, \frac {h-1} {h+1} \bigg) \bigg\}.
\end{equation}
\end{proposition}

For the external field $Q$ generated by a point charge outside of the sphere $\S^2$, we can, in fact, specify when the extremal support $S_Q$ is a \textit{proper} subset of the sphere $\S^2$. Indeed, the answer is given in terms of solution to Gonchar's problem, as explained in \cite{bds3}, \cite{bdsw}. It states that if a point charge, located on the positive polar semi-axis, is too far from the surface of the sphere $\S^2$, or the magnitude $q$ of the point charge is too small, the electrostatic field, created by this point charge, is too weak to force the equilibrium distribution from occupying the whole surface of the sphere $\S^2$. More specifically, in order for the extremal support $S_Q$ to be a spherical cap $C_{S,\alpha}$, the height of the point charge of magnitude $q$ must satisfy $h<h_{+}$, where $h_{+}$ is defined as the unique (real) zero in $(1,\infty)$ of the cubic equation
\begin{equation}\label{sgp}
\frac {1} {q} = \frac {h+1} {(h-1)^2} - \frac {1} {h}.
\end{equation}
For instance, in the case of the positive unit point charge $q=1$, the critical height $h_{+}$ is precisely $1+\rho_+$, where $\rho_+$ is the golden ratio $(1+\sqrt{5})/2  \approx 1.6180339887$.

If a positive point charge of magnitude $q$ is placed inside the sphere $\S^2$, then in order that the extremal support $S_Q$ to be a spherical cap $C_{S,\alpha}$, the height of the point charge of magnitude $q$ must satisfy $h>h_{-}$, where $h_{-}$ is defined as the unique (real) zero in $(0,1)$ of the cubic equation (see \cite{bds3})
\begin{equation}\label{sgp1}
\frac {1} {q} = \frac {h+1} {(1-h)^2} - 1.
\end{equation}
We next write down a transcendental equation, whose solution gives the angle $\alpha$, defining the support $C_{S,\alpha}$ of the extremal measure $\mu_Q$.
\begin{theorem}\label{transceqsupp}
For the external field $Q$ of $(\ref{ptch})$, with $h$ chosen such that $h<h_{+}$ or $h>h_{-}$, the support $S_Q$ is a spherical cap $C_{S,\alpha}$, centered at the South Pole, with $\alpha=\alpha_0\in(0,\pi)$. The angle $\alpha_0$ is defined as a unique solution of the transcendental equation
\begin{empheq}{align}\label{alsupp}
\frac {\pi} {\pi-\alpha+\sin \alpha}  \bigg\{1 + \frac {q(h+1)}{2h} \bigg(1 - \frac {\alpha} {\pi} \bigg)   - \frac {q(h-1)}{\pi h} \tan^{-1}\bigg(  \cot & \bigg(  \frac{\alpha}{2}\bigg)  \, \frac {h-1} {h+1} \bigg) \bigg\} \\ \nonumber
 & =  q \, \frac {h+1} {h^2+1-2h\cos \alpha}.
\end{empheq}
\end{theorem}
Now, as the support $C_{S,\alpha}$ is found via equation $(\ref{alsupp})$, we can recover the equilibrium measure $\mu_Q$ from integral equation $(\ref{ieec})$. Using Theorem \ref{theo3}, by a direct calculation we obtain the following result.
\begin{theorem}\label{eqdgct}
Let $\alpha_0$ be a unique solution to $(\ref{alsupp})$. Set 
\begin{empheq}{align}\label{eqdgcf}
F(\varphi) & = - \frac {q(h+1)} {2\pi^2} \bigg\{ \frac{1}{1+h^2-2h\cos\varphi} \sqrt{ \frac{1-\cos \alpha_0} {\cos \alpha_0-\cos\varphi}} \\ \nonumber
& + \frac {h-1}{(1+h^2-2h\cos\varphi)^{3/2}}  \tan^{-1} \bigg( \frac{h-1}{\sqrt {1+h^2-2h\cos\varphi}}\sqrt{ \frac{\cos \alpha_0-\cos\varphi} {1-\cos \alpha_0)}} \bigg) \bigg\},
\end{empheq}
\begin{equation}\label{FQgc}
F_Q = \frac {\pi} {\pi-\alpha_0+\sin \alpha_0} \bigg\{ 1 + \frac {q(h+1)}{2h} \bigg(1 - \frac {\alpha_0} {\pi} \bigg) - \frac {q(h-1)}{\pi h} \tan^{-1}\bigg( \frac {h-1} {h+1} \cot\bigg(\frac{\alpha_0}{2}\bigg) \bigg) \bigg\}.
\end{equation}
The density $f$ of the equilibrium measure $\mu_Q$ is given by
\begin{equation}\label{eqds1}
f(\varphi) = \frac {F_Q} {4\pi} \left\{ 1 + \frac {2} {\pi} \left[  \sqrt{ \frac {1- \cos \alpha_0} {\cos \alpha_0 - \cos\varphi} } -\tan^{-1} \sqrt{\frac {1 - \cos \alpha_0} {\cos \alpha_0 - \cos\varphi} }   \right]    \right\} + F(\varphi),
\end{equation}
where $\alpha_0 \leq \varphi \leq \pi $.
\end{theorem}

From Theorem \ref{transceqsupp} it follows that in the case when the point charge is placed at the North Pole, equation $(\ref{alsupp})$ takes especially simple form. From $(\ref{eqdgcf})$ it transpires that in this case the equilibrium density will have a simple structure as well.
\begin{theorem}\label{chargenpole}
For the external field $Q$ of $(\ref{ptch})$, with the point charge placed at the North Pole, the support $S_Q$ is a spherical cap $C_{S,\alpha}$, centered at the South Pole, with $\alpha=\alpha_0\in(0,\pi)$. The angle $\alpha_0$ is defined as a unique solution of the transcendental equation
\begin{empheq}{align}\label{alsuppnp}
\pi \sec(\alpha) = \pi +q(\pi-\alpha) + q \tan \alpha.
\end{empheq}
The equilibrium density is given in this case by
\begin{empheq}{align}\label{edh1}
f(\varphi) = & \frac {\pi+q(\pi-\alpha_0)} {4\pi (\sin \alpha_0+\pi - \alpha_0)} \left\{ 1 + \frac {2} {\pi} \left[  \sqrt{ \frac {1- \cos \alpha_0} {\cos \alpha_0 - \cos\varphi} } -\tan^{-1} \sqrt{\frac {1 - \cos \alpha_0} {\cos \alpha_0 - \cos\varphi} }   \right]    \right\} \\ \nonumber
- &  \frac {q} {2\pi^2}  \frac{1} {1-\cos\varphi} \sqrt{\frac{1-\cos \alpha_0}{\cos \alpha_0-\cos\varphi}}, \quad \alpha_0\leq \varphi \leq \pi.
\end{empheq}
\end{theorem}

The second application of our results is concerned with the quadratic external fields. In particular, we consider a quadratic external field
\begin{equation}\label{polyextf}
Q(x)=Q(x_3) = a x_3^2 + b x_3 +c, \quad x\in\S^2,
\end{equation}
where coefficients $a,b,c$ chosen such that $4a^2 < b^2 \leq 4ac$, with $a,b>0$. The choice of the coefficients implies that the external field $Q$ is a convex increasing function on $[-1,1]$, while its invariance with respect to rotations about the polar axis $x_3$ is obvious. It then follows from Proposition \ref{scssp} that the support of the extremal measure is a cap $C_{S,\alpha}$, centered at the South Pole. Calculating $\cF$-functional for such $Q$, we obtain the following statement.
\begin{proposition}\label{ffuncpolyextf}
The $\cF$-functional for the spherical cap $C_{S,\alpha}$, in the case of the quadratic external field $Q$, as in $(\ref{polyextf})$, is given by
\begin{empheq}{align}\label{msfexpquadpoly}
\cF(C_{S,\alpha})  & =  \frac {1} {36(\pi-\alpha+\sin \alpha)}  \bigg\{  \tan \left( \frac{\alpha}{2}\right) \bigg [32\, a \cos^3 \alpha  + 4(2a+9b) \cos^2 \alpha \\  \nonumber
& + 4(9c-5a) \cos \alpha + 4 a - 36 b + 36 c \bigg]  + 12 (a+3c) (\pi-\alpha) + 36 \pi \bigg\}.
\end{empheq}
\end{proposition}

\noindent The support $C_{S,\alpha}$ of the extremal measure $\mu_Q$ is determined via a transcendental equation, as follows.
\begin{theorem}\label{transceqsupppolyextf}
For the quadratic external field $Q$ of $(\ref{polyextf})$, the support $S_Q$ is a spherical cap $C_{S,\alpha}$, centered at the South Pole, with $\alpha=\alpha_0\in(0,\pi)$. The angle $\alpha_0$ is defined as a unique solution of the transcendental equation
\begin{empheq}{align}\label{alsupppoly}
 8a\cos^3 \alpha (2\sin\alpha+3(\pi-\alpha)) & + \cos^2\alpha((2a+9b)\sin\alpha - 6(2a-3b)(\pi-\alpha)) \\ \nonumber
& + \frac{1}{2}\sin(2\alpha) (9b-22a) + 3(2a-3b)(\pi-\alpha)+9\pi \\ \nonumber 
& = 9\cos\alpha (\pi + (2a+b)(\pi-\alpha)) - 2\sin\alpha(2a-9b).
\end{empheq}
\end{theorem}

\noindent Having the support $C_{S,\alpha}$ now determined, we compute the density of the equilibrium measure $\mu_Q$ by invoking Theorem \ref{theo3}.
\begin{theorem}\label{eqdtheopoly}
Let $\alpha_0$ be a unique solution to $(\ref{alsupppoly})$. Set 
\begin{empheq}{align}\label{eqdpoly}
F(\varphi)  = & \frac {1}{36\,\pi^2} \bigg\{  \sqrt{1-\cos\alpha_0}\, \sqrt{\cos\alpha_0-\cos\varphi}\, (20a\cos\alpha_0 + 60a\cos\varphi + 10a + 27b) \\ \nonumber
 & - \sqrt{\frac{1-\cos\alpha_0}{\cos\alpha_0-\cos\varphi}}\, (8a\cos^2\alpha_0 + 10a\cos\alpha_0 \cos\varphi + (4a+9b) \cos\alpha_0 \\ \nonumber
		& + (20a + 27b) \cos\varphi + 15a\cos(2\varphi) + 9a + 18b + 18c) \\ \nonumber
		& - 6 \tan^{-1} \sqrt{\frac{\cos\alpha_0-\cos\varphi}{1-\cos\alpha_0}}\, (15a\cos^2\varphi + 9b\cos\varphi - 4a + 3c) \bigg\},
\end{empheq}
\begin{empheq}{align}\label{FQpoly}
F_Q & =  \frac {1} {36(\pi-\alpha_0+\sin \alpha_0)}  \bigg\{  \tan \left( \frac{\alpha_0}{2}\right) \bigg [32\, a \cos^3 \alpha_0  + 4(2a+9b) \cos^2 \alpha_0 \\  \nonumber
& + 4(9c-5a) \cos \alpha_0 + 4 a - 36 b + 36 c \bigg]  + 12 (a+3c) (\pi-\alpha_0) + 36 \pi \bigg\}.
\end{empheq}
The density $f$ of the equilibrium measure $\mu_Q$ is given by
\begin{equation}\label{eqdspoly}
f(\varphi) = \frac {F_Q} {4\pi} \left\{ 1 + \frac {2} {\pi} \left[  \sqrt{ \frac {1- \cos \alpha_0} {\cos \alpha_0 - \cos\varphi} } -\tan^{-1} \sqrt{\frac {1 - \cos \alpha_0} {\cos \alpha_0 - \cos\varphi} }   \right]    \right\} + F(\varphi),
\end{equation}
where $\alpha_0 \leq \varphi \leq \pi $.
\end{theorem}


\section{Proofs}
\noindent{\bf Proof of Theorem \ref{theo1}.} The strategy of the proof is to show that the equilibrium potential $U^{\mu_Q}$ is Lipschitz continuous on $C_{S,\alpha}$. 
If that is established, it will imply that the normal derivatives of $U^{\mu_Q}$ exist a.e. on  $C_{S,\alpha}$. Hence the 
measure $\mu_Q$ can be recovered from its potential by the formula
\begin{equation*}
d\mu_Q = -\frac{1}{4\pi}\left(\frac{\partial U^{\mu_Q}}{\partial n_+}+\frac{\partial U^{\mu_Q}}{\partial n_-}\right)\,dS:=f(\varphi)\, \sin\varphi\,d\theta\,d\varphi,
\end{equation*}
where $dS= \sin\varphi\,d\theta\,d\varphi$ is the Lebesgue surface measure on $\supp(\mu_Q)$, $n_{+}$  and 
$n_{-}$ are the inner and the outer normals to the cap $C_{S,\alpha}$, see \cite[pp. 164--165]{land} and \cite[p. 164]{kel}. It is clear that 
the normal derivatives of $U^{\mu_Q}$ are bounded a.e. by the Lipschitz constant, and hence we 
obtain $f\in L_{\text{loc}}^\infty([\alpha,\pi])$. 

Our first step is to construct an extension 
of the external field $Q$ to $\S^2$ in such a way that the extremal measures for the cap $C_{S,\alpha}$ and 
the sphere $\S^2$ are the same. In doing so, we will follow an idea of Pritsker \cite{prit}. Recall that the external field $Q$ is a $C^2$ function on an open neighborhood $N$ of $S_Q$ in $\S^2$. We can adjust $Q$ in such a way that 
for the new external field $\widetilde{Q}$ one has
\[
U^{\mu_Q}(x)+\widetilde{Q}(x) = F_Q,  \quad x \in S_Q,
\]
\[
U^{\mu_Q}(x)+\widetilde{Q}(x) \geq F_Q,  \quad x \in \S^2,
\]
and also $\widetilde{Q}\in C^2(\S^2)$.
To show that it is indeed possible, we first remark that the external field $Q$ is rotationally symmetric about the polar axis. 
Therefore, $Q$ is a function of the polar angle $\varphi$ only, that is $Q=Q(\varphi)$. This symmetry 
is also inherited by potential, hence on the sphere $\S^2$ we have $U^{\mu_Q}(x)=U^{\mu_Q}(\varphi)$, $x=(r,\theta,\varphi)\in\S^2$. Next, note that $N= \{ (r,\theta,\varphi): r=1,\ \ \gamma <  \varphi \leq \pi,\ \  0\leq \theta \leq 2\pi  \}$, with some $\gamma\in(0,\alpha)$. Pick a number $\epsilon$ such that $\gamma < \epsilon < \alpha$. We define a new external field $\widetilde{Q}$ as follows: set $\widetilde{Q}(\varphi)=Q(\varphi)$, for $\epsilon < \varphi \leq \pi$, while on  $[0, \epsilon]$ we tweak $Q$ to $\widetilde{Q}$ in such a way that 
\[
U^{\mu_Q}(\varphi)+\widetilde{Q}(\varphi)\geq F_Q,
\]
and $\widetilde{Q}\in C^2(\S^2)$. Applying Theorem 4.2.14 from \cite{bhs}, we infer that $\mu_{\widetilde{Q}}=\mu_Q$ and $F_{\widetilde Q}=F_Q$. 

Let $u$ and $v$ denote the equilibrium potentials for the minimum energy problem on $C_{S,\alpha}$ and $\S^2$, 
respectively. Since the equilibrium measure is the same for those two sets, it immediately follows that $u=v$. 
Now observe that the spherical cap $C_{S,\alpha}$ is a part of the sphere $\S^2$, which is a closed smooth surface. Thus 
we can invoke the result of G\"otz \cite{go} to conclude that $v$ 
is $C^{2}$ in an open neighborhood $\mathcal{U}$ of $\S^2$. Now recalling that $u=U^{\mu_Q}$, we 
deduce that $U^{\mu_Q}$ is Lipschitz continuous on $C_{S,\alpha}$.

\qed

\noindent{\bf Proof of Theorem \ref{theo2}.} Assume that $S_Q = C_{N,\alpha}$. From Theorem \ref{theo1} we know 
that $d\mu_Q = f(\varphi)\, \sin\varphi\,d\theta\,d\varphi $, where $f\in L_{\text{loc}}^\infty([0,\alpha])$. 

In what follows, we will need an expression for the distance between a point on a sphere $\S^2$ and another point in the space $\R^3$, which is not on the surface of the sphere $\S^2$. Let
\begin{align*}
& x=(r\sin \varphi_1\cos \theta_1,r\sin \varphi_1\sin \theta_1,r\cos \varphi_1)\in\R^3,\\
& y=(\sin \varphi_2\cos \theta_2,\sin \varphi_2\sin \theta_2,\cos \varphi_2)\in\S^2,
\end{align*}
be such two points, written in spherical coordinates. Then, for the distance $|x-y|$, one has
\begin{align*}
 |x-y|^2 & = r^2 \sin^2 \varphi_1 \cos^2 \theta_1 + \sin^2 \varphi_2 \cos^2 \theta_2-2r\sin \varphi_1 \sin \varphi_2 \cos \theta_1\cos \theta_2 \\
 		& + r^2 \sin^2 \varphi_1 \sin^2 \theta_1 + \sin^2 \varphi_2 \sin^2 \theta_2-2r\sin \varphi_1 \sin \varphi_2 \sin \theta_1\sin \theta_2 \\
		& + r^2\cos^2 \varphi_1 + \cos^2 \varphi_2 - 2r\cos \varphi_1\cos \varphi_2 \\
		& = r^2 + 1 - 2r(\cos \varphi_1\cos \varphi_2+\sin \varphi_1\sin \varphi_2\cos(\theta_1-\theta_2)) \\
		& = r^2 + 1 - 2r\gamma,
\end{align*}
where $\gamma = \cos \varphi_1\cos \varphi_2+\sin \varphi_1\sin \varphi_2\cos(\theta_1-\theta_2)$.

\noindent The potential $U^{\mu_Q}$ of 
$\mu_Q$ has the form
\begin{equation}\label{pinspc}
U^{\mu_Q}(r,\theta,\varphi) = \int_0^\alpha \int_0^{2\pi} \frac {f(\xi)\sin \xi \,d\eta \, d\xi } {\sqrt{1+r^2-2\gamma r } },
\end{equation}
where $\gamma = \cos\varphi\cos \xi+\sin\varphi\sin \xi\cos(\theta-\eta)$.

On the sphere $\S^2$, we have $r=1$. Hence, from $(\ref{pinspc})$ it is clear that the potential  $U^{\mu_Q}$  on the sphere $\S^2$ is 
\[
U^{\mu_Q}(\theta,\varphi) = \int_0^\alpha \int_0^{2\pi} \frac {f(\xi)\sin \xi \,d\eta \, d\xi } {\sqrt{2-2\gamma}}.
\]
Equation $(\ref{ieec})$ now takes the form
\begin{equation}\label{ie1}
\int_0^\alpha \int_0^{2\pi} \frac {f(\xi)\sin \xi \,d\eta \, d\xi } {\sqrt{2-2\gamma}} = F_Q - Q(\varphi), \quad 0 \leq \varphi \leq \alpha,
\end{equation}
where $\gamma = \cos\varphi\cos \xi+\sin\varphi\sin \xi\cos(\theta-\eta)$.

Equation $(\ref{ie1})$ was first considered and solved in \cite{col}, in connection with the problem of electrified 
spherical cap. The solution of this equation is based on ideas of Copson \cite{cop}, who considered a similar problem
for the disk in the plane. We present the derivation of solution to $(\ref{ie1})$, which is more tailored to our needs, here.

The first step is to transform a somewhat cumbersome term $2-2\gamma$. First, note that
\begin{align*}
\sin\varphi\sin \xi & =  2\sin\left( \frac {\varphi} {2} \right)  \cos\bigg( \frac {\xi} {2} \bigg)  \, 2 \sin\bigg( \frac {\xi} {2} \bigg) \cos\left( \frac {\varphi} {2} \right).
\end{align*}
Next, using the half-angle identities, it can be shown that
\begin{align*}
1 - \cos\varphi\cos \xi   = 2\left(\sin^2 \left( \frac {\varphi} {2} \right)\cos^2\bigg( \frac {\xi} {2} \bigg)  + \sin^2\bigg( \frac {\xi} {2} \bigg)  \cos^2 \left( \frac {\varphi} {2} \right) \right).
\end{align*}
\noindent Therefore, we see that
\begin{align*}
2 - 2\gamma & = 4 \left(\sin^2 \left( \frac {\varphi} {2} \right)\cos^2\bigg( \frac {\xi} {2} \bigg)  + \sin^2\bigg( \frac {\xi} {2} \bigg)  \cos^2 \left( \frac {\varphi} {2} \right) \right) \\
			& - 2\cdot \, 2\sin\left( \frac {\varphi} {2} \right)  \cos\bigg( \frac {\xi} {2} \bigg)  \, 2 \sin\bigg( \frac {\xi} {2} \bigg) \cos\left( \frac {\varphi} {2} \right)\, \cos(\theta-\eta) \\
			& = a^2 + b^2 - 2ab \cos(\theta-\eta),
\end{align*}
where
\[
a:= 2\sin\left( \frac {\varphi} {2} \right)  \cos\bigg( \frac {\xi} {2} \bigg),
\]
\[
b:= 2 \sin\bigg( \frac {\xi} {2} \bigg) \cos\left( \frac {\varphi} {2} \right).
\]
Note that since $0<\alpha\leq\pi$, it is clear that $a>0$ and $b>0$ for all and $0\leq \varphi\leq \alpha$ and $0\leq \xi\leq\alpha$. We need a special case of Copson's lemma \cite[p. 16]{cop}.
\begin{lemma}\label{coplem}
If $a$ and $b$ are positive numbers, then
\begin{equation}\label{copf}
\int_0^{2\pi} \frac {d\eta} {\sqrt{a^2+b^2-2ab\cos\eta}} = 4 \int_0^{\min{(a,b)}} \frac {dt} {\sqrt{a^2-t^2} \sqrt{b^2-t^2}}.
\end{equation}
\end{lemma}
\noindent In view of Lemma \ref{coplem}, the inner integral in left hand side of equation $(\ref{ie1})$ is transformed as
\begin{align*}
 \int_0^{2\pi} \frac {d\eta } {\sqrt{2-2\gamma}} & = \int_0^{2\pi} \frac {d\eta} {\sqrt{a^2 + b^2 - 2ab \cos(\theta-\eta)}} \\
									& = 4 \int_0^{\min{(a,b)}} \frac {dt} {\sqrt{a^2-t^2} \sqrt{b^2-t^2}}.
\end{align*}
One can easily check that $a<b$ for $\varphi<\xi$, while for $\varphi>\xi$, we have $a>b$. Splitting the interval of integration of the outer integral in the 
left hand side of equation $(\ref{ie1})$, we rewrite equation $(\ref{ie1})$ as follows,
\begin{align*}
 \int_0^\varphi f(\xi) \sin \xi\, d\xi  \int_0^b \frac {dt} {\sqrt{a^2-t^2} \sqrt{b^2-t^2}} + & \int_\varphi^\alpha f(\xi) \sin \xi\, d\xi  \int_0^a \frac {dt} {\sqrt{a^2-t^2} \sqrt{b^2-t^2}} \\ = & \frac {1} {4} \, (F_Q - Q(\varphi)), \quad 0 \leq \varphi \leq \alpha.
\end{align*}
By virtue of the substitution
\[
t = 2 \cos\bigg(\frac {\varphi} {2}\bigg) \cos\bigg( \frac {\xi} {2}\bigg) \tan\bigg( \frac {\zeta} {2}\bigg),
\]
we can further obtain
\begin{empheq}{align}\label{ie2}
& \frac{1} {2} \int_0^\varphi f(\xi) \sin \xi \, d\xi  \int_0^\xi \frac {d\zeta} {\sqrt{\cos \zeta - \cos\varphi} \sqrt{\cos \zeta - \cos \xi}} \\\nonumber
		 + & \frac{1} {2} \int_\varphi^\alpha f(\xi) \sin \xi  \, d\xi  \int_0^\varphi \frac {d\zeta} {\sqrt{\cos \zeta - \cos\varphi} \sqrt{\cos \zeta - \cos \xi}}  = \frac {1} {4} \, (F_Q - Q(\varphi)), \quad 0 \leq \varphi \leq \alpha.
\end{empheq}
Inverting the order of integration in the first integral in the left hand side of $(\ref{ie2})$, as
\begin{align*}
& \int_0^\varphi f(\xi) \sin \xi  \, d\xi  \int_0^\xi \frac {d\zeta} {\sqrt{\cos \zeta - \cos\varphi} \sqrt{\cos \zeta - \cos \xi}} \\
		 = & \int_0^\varphi \frac {d\zeta} {\sqrt{\cos \zeta - \cos\varphi}} \int_\zeta^\varphi \frac {f(\xi)  \sin \xi  \, d\xi} {\sqrt{\cos \zeta - \cos \xi}},    
\end{align*}
we recast equation $(\ref{ie2})$ into
\begin{equation}\label{ie3}
 \int_0^\varphi \frac {d\zeta} {\sqrt{\cos \zeta - \cos\varphi}} \int_\zeta^\alpha  \frac{1} {2}  \frac {f(\xi)  \sin \xi \, d\xi} {\sqrt{\cos \zeta - \cos \xi}}  =  \frac {1} {4} \, (F_Q - Q(\varphi)), \quad 0 \leq \varphi \leq \alpha.
\end{equation}
Letting
\begin{equation}\label{ie4}
S(\zeta) =  \int_\zeta^\alpha  \frac{1} {2}  \frac {f(\xi)  \sin \xi  \, d\xi} {\sqrt{\cos \zeta - \cos \xi }}\, \quad 0 \leq \zeta \leq \alpha,
\end{equation}
we see that equation $(\ref{ie3})$ becomes
\begin{equation}\label{ie5}
 \int_0^\varphi \frac {S(\zeta)\, d\zeta} {\sqrt{\cos \zeta - \cos\varphi}} = \frac {1} {4} \, (F_Q - Q(\varphi)), \quad 0 \leq \varphi \leq \alpha.
\end{equation}
Equation $(\ref{ie5})$ is an Abel type integral equation. Since $Q\in C^2$, the solution of $(\ref{ie5})$ is \cite[p. 50, \# 23]{polman}
\begin{equation}
S(\zeta) = \frac {1} {4\pi} \frac {d} {d\zeta} \int_0^\zeta \frac {(F_Q - Q(\varphi)) \sin\varphi \, d\varphi} {\sqrt{\cos\varphi - \cos \zeta}}, \quad 0 \leq \zeta \leq \alpha.
\end{equation}
Observing that equation $(\ref{ie4})$ is again an Abel type integral equation with respect to $f(\xi) \sin \xi$, we solve it and obtain
\begin{equation}\label{ie6}
f(\xi) = - \frac {2} {\pi} \frac {1} {\sin \xi} \frac {d} {d\xi} \int_\xi^\alpha \frac {S(\zeta) \sin \zeta\, d\zeta} {\sqrt{\cos \xi - \cos \zeta}}, \quad 0  \leq \xi \leq \alpha.
\end{equation}
Denote
\begin{equation}\label{ie7}
g(\zeta) = \frac {1} {4\pi}\, \frac {d} {d\zeta} \int_0^\zeta \frac{\sin \xi Q(\xi) \,d\xi} {\sqrt{\cos \xi - \cos \zeta}},
\end{equation}
and let
\begin{equation}\label{ie8}
F(\xi) = \frac {2} {\pi}\frac {1} {\sin \xi}  \frac {d} {d\xi} \int_\xi^\alpha \frac {g(\zeta) \sin \zeta \, d\zeta} {\sqrt{\cos \xi - \cos \zeta}}, \quad 0  < \xi < \alpha.
\end{equation}
In view of $(\ref{ie7})$ and $(\ref{ie8})$, expression $(\ref{ie6})$ takes the form
\begin{equation}\label{ie9}
f(\xi)  =  - \frac {F_Q} {2\pi^2} \frac {1} {\sin \xi}  \frac {d} {d\xi} \int_\xi^\alpha   \frac {1} {\sqrt{\cos \xi - \cos \zeta}} \bigg\{ \frac {d} {d\zeta} \int_0^\zeta \frac{\sin \xi\,d\xi} {\sqrt{\cos \xi - \cos \zeta}} \bigg\}\, d\zeta  + F(\xi), \quad 0  \leq \xi \leq \alpha.
\end{equation}
Executing the sequence of integrations and differentiations on the right hand side of $(\ref{ie9})$, we arrive to
\begin{equation}\label{dfchk}
f(\xi) = \frac {F_Q}{4\pi} \bigg\{ 1 + \frac{2}{\pi} \bigg[  \sqrt{ \frac {1+\cos \alpha} {\cos \xi - \cos \alpha} } - \tan^{-1} \sqrt{ \frac {1+\cos \alpha} {\cos \xi - \cos \alpha}}  \bigg] \bigg\} + F(\xi),\end{equation}
which is the desired expression $(\ref{eqd})$.

We now show how to find the Robin constant $F_Q$. Recall that $\mu_Q$ has mass $1$, so that
\begin{equation}\label{mmo}
1=\int d\mu_Q = \int_0^\alpha \int_0^{2\pi} f(\varphi)\, \sin\varphi\,d\theta\, d\varphi.
\end{equation}
Inserting expression $(\ref{dfchk})$ into $(\ref{mmo})$, we obtain
\begin{empheq}{align}\label{constfq}
1 & = 2\pi \int_0^\alpha F(\varphi) \, \sin\varphi \, d\varphi \\ \nonumber
	& + \frac{F_Q}{2}  \int_0^\alpha \bigg\{ 1 + \frac{2}{\pi} \bigg[  \sqrt{ \frac {1+\cos \alpha} {\cos\varphi - \cos \alpha} } - \tan^{-1} \sqrt{ \frac {1+\cos \alpha} {\cos\varphi - \cos \alpha} }  \bigg] \bigg\} \sin\varphi \, d\varphi.
\end{empheq}
It is easy to see that
\[
\int_0^\alpha \sqrt{ \frac {1+\cos \alpha} {\cos\varphi - \cos \alpha} } \sin\varphi \, d\varphi = 2 \sin \alpha.
\]
Making the change of variables $t=(\cos\varphi-\cos \alpha)/(1+\cos \alpha)$, and then integrating by parts, one can show that
\[
\int_0^\alpha \tan^{-1} \sqrt{ \frac {1+\cos \alpha} {\cos\varphi - \cos \alpha} } \sin\varphi \, d\varphi = \frac {\pi (1-\cos \alpha)} {2} + \sin \alpha - \alpha.
\]
Using the last two facts, we see that 
\begin{align*}
 \int_0^\alpha \bigg\{ 1 + \frac{2}{\pi} \bigg[  \sqrt{ \frac {1+\cos \alpha} {\cos\varphi - \cos \alpha} }  & -  \tan^{-1} \sqrt{ \frac {1+\cos \alpha} {\cos\varphi - \cos \alpha} }  \bigg] \bigg\}  \sin\varphi  \, d\varphi   \\
& = 1 - \cos \alpha + \frac {4} {\pi} \sin \alpha - (1 - \cos \alpha) - \frac {2} {\pi} \sin \alpha + \frac {2} {\pi} \alpha \\
& =  \frac{2}{\pi}(\sin \alpha + \alpha).
\end{align*}
Therefore, 
\[
1 = 2\pi \int_0^\alpha F(\varphi) \, \sin\varphi \, d\varphi + \frac{F_Q}{\pi}(\sin \alpha + \alpha),
\]
and thus
\[
F_Q = \frac {\pi} {\alpha+\sin \alpha} \bigg\{ 1 - 2\pi \int_0^\alpha F(\varphi) \sin\varphi \, d\varphi \bigg\}.
\] 
This is formula $(\ref{dcf})$, as desired.

\qed


\noindent{\bf Proof of Theorem \ref{theo3}.} If $S_Q=C_{S,\alpha}$, equation $(\ref{ieec})$ assumes the form
\begin{equation}\label{ie21}
\int_\alpha^\pi \int_0^{2\pi} \frac {f(\xi)\sin \xi \,d\eta \, d\xi } {\sqrt{2-2\gamma}} = F_Q - Q(\varphi), \quad \alpha \leq \varphi \leq \pi ,
\end{equation}
where $\gamma = \cos\varphi\cos \xi+\sin\varphi\sin \xi\cos(\theta-\eta)$.
Via the change of variables $\widetilde{\varphi} = \pi - \varphi$, we transform $(\ref{ie21})$ into 
\begin{equation}\label{ie22}
\int_0^\beta \int_0^{2\pi} \frac {f_0(\xi)\sin \xi \,d\eta \, d\xi } {\sqrt{2-2\widetilde{\gamma}}} = F_Q - Q_0(\varphi), \quad 0 \leq \widetilde{\varphi} \leq \beta,
\end{equation}
with $\beta=\pi-\alpha$, $f_0(\widetilde{\varphi})=f(\pi-\widetilde{\varphi})$ and $Q_0(\widetilde{\varphi})=Q(\pi-\widetilde{\varphi})$, and $\widetilde{\gamma} = \cos \widetilde{\varphi} \cos \xi+\sin \widetilde{\varphi} \sin \xi\cos(\theta-\eta)$. The integral equation $(\ref{ie22})$ is of the form $(\ref{ie1})$. Hence Theorem \ref{theo2} applies, and we obtain
\begin{equation*}
F(\widetilde{\varphi}) = \frac {2} {\pi} \frac {1} {\sin \widetilde{\varphi} } \, \frac {d} {d\widetilde{\varphi}} \int_{\widetilde{\varphi}}^\beta 
\frac{\sin t g(t)\,dt} {\sqrt{\cos \widetilde{\varphi}  - \cos t }},
\end{equation*}
where
\begin{equation*}
g(t) = \frac {1} {4\pi}\, \frac {d} {dt} \int_0^t \frac{\sin \xi Q_0(\xi) \,d\xi} {\sqrt{\cos \xi - \cos t}}.
\end{equation*}
The density $f$ of the equilibrium measure $\mu_Q$ in this case is given by
\begin{equation*}
f(\widetilde{\varphi}) = \frac {F_Q} {4\pi} \left\{ 1 + \frac {2} {\pi} \left[  \sqrt{ \frac {1+\cos \beta} {\cos \widetilde{\varphi}  - \cos \beta} } -\tan^{-1} \sqrt{ \frac {1+\cos \beta} {\cos \widetilde{\varphi}  - \cos \beta} }   \right]    \right\} + F(\widetilde{\varphi}),
\end{equation*}
where $0 \leq \widetilde{\varphi} \leq \beta $. Going back to the $\varphi$ variable via $\varphi = \pi - \widetilde{\varphi}$, after some algebra, we find
\begin{equation*}
g(t) = \frac {1} {4\pi}\, \frac {d} {dt} \int_t^\pi \frac{\sin \xi Q(\xi) \,d\xi} {\sqrt{\cos t - \cos \xi}},
\end{equation*}
and thus
\begin{equation*}
F(\varphi) = \frac {2} {\pi} \frac{1} {\sin\varphi} \, \frac {d} {d\varphi} \int_\alpha^{\varphi} 
\frac{\sin t \, g(t)\,dt} {\sqrt{\cos t - \cos\varphi }}.
\end{equation*}
For the Robin constant $F_Q$, in this case we similarly obtain the expression
\begin{equation*}
F_Q = \frac {\pi} {\sin \alpha + \pi - \alpha} \left\{ 1 - 2\pi \int_\alpha^\pi F(\varphi) \sin\varphi\,d\varphi \right\}.
\end{equation*}

\qed


\noindent{\bf Proof of Proposition \ref{msfptch}.} The external field $Q$, generated by a positive point charge
of magnitude $q$, located on the polar axis at the distance $h$ from the origin, has the following expression in
terms of spherical coordinates,
\begin{equation}\label{extfptch}
Q(\varphi) = \frac {q} {\sqrt{1+h^2-2h\cos\varphi}}.
\end{equation}
We now insert expression $(\ref{extfptch})$ into $(\ref{mhsf})$, and evaluate the integrals on the right hand side 
of $(\ref{mhsf})$. The first integral on the right hand side of $(\ref{mhsf})$ evaluates as follows,
\begin{empheq}{align}\label{i1}
 \int_\alpha^\pi Q(\varphi) \sin\varphi\, d\varphi & = q  \int_\alpha^\pi \frac{ \sin\varphi\, d\varphi }  {\sqrt{1+h^2-2h\cos\varphi}} \\ \nonumber
 									& = \frac {q} {h} \bigg\{ (1+h) - \sqrt{1+h^2-2h\cos \alpha } \bigg\}.
\end{empheq}
Next, we deal with the second integral on the right hand side of $(\ref{mhsf})$. We have
\begin{empheq}{align}\label{i2}
 \int_\alpha^\pi Q(\varphi) \sqrt{ \frac {1-\cos \alpha} {\cos \alpha - \cos\varphi} }  \sin\varphi\,  & d\varphi   = q \sqrt{1-\cos \alpha} \int_\alpha^\pi  \frac{ \sin\varphi\, d\varphi }  {\sqrt{1+h^2-2h\cos\varphi}\, \sqrt{\cos \alpha - \cos\varphi}} \\ \nonumber
& = \frac{2q}{\sqrt{h}}\, \sin\bigg(\frac{\alpha}{2}\bigg) \,  \log\bigg( \frac {1+h+2\sqrt{h}\cos(\alpha/2)} {\sqrt{1+h^2-2h\cos \alpha}} \bigg).
\end{empheq}
We now turn to the last remaining integral on the right hand side of $(\ref{mhsf})$. Namely,
\begin{align*}
 \int_\alpha^\pi  Q(\varphi) \tan^{-1}\sqrt{ \frac {1-\cos \alpha} {\cos \alpha - \cos\varphi}}  \sin\varphi\, d\varphi & = q \int_\alpha^\pi  \tan^{-1} \sqrt{ \frac {1-\cos \alpha} {\cos \alpha - \cos\varphi} } \, \frac{ \sin\varphi\, d\varphi} {\sqrt{1+h^2-2h\cos\varphi}} \\
& = \frac {q} {\sqrt{2h}} \int_{-1}^a \frac {1} {\sqrt{c-u}}  \tan^{-1} \sqrt{ \frac {1-a} {a - u} } \, du,
\end{align*}
where we set $u=\cos\varphi$, $a=\cos \alpha$, and $c = (h^2+1)/2h$. Integrating by parts, and then using the substitution $a-u=t^2$, it can be showed that
\begin{align*}
&  \int_{-1}^a \frac {1} {\sqrt{c-u}}  \tan^{-1} \sqrt{ \frac {1-a} {a - u}} \, du  = - \pi \sqrt{c-a} + 2 \sqrt{c+1} \tan^{-1}\bigg( \frac {1-a} {1 + a} \bigg)^{1/2}  \\
 	& + 2 \sqrt{1-a} \bigg\{ \sqrt{\frac{c-1}{1-a}} \tan^{-1}  \sqrt{ \frac {(1+a)(c-1)} {(1-a)(c+1)} } + \log (\sqrt{c+1}+\sqrt{a+1}) - \frac {1} {2} \log(c-a) \bigg\}.
\end{align*}
After some algebra, we can conclude that
\begin{align*}
 \int_\alpha^\pi & Q(\varphi) \tan^{-1} \sqrt{ \frac {1-\cos \alpha} {\cos \alpha - \cos\varphi} }  \sin\varphi\, d\varphi \\
 & = \frac {q} {2h} \bigg\{ (h+1)\alpha + 2(h-1) \tan^{-1}\bigg( \cot\bigg(\frac{\alpha}{2}\bigg) \, \frac {h-1} {h+1} \bigg) - \pi \sqrt{1+h^2 - 2h\cos \alpha} \bigg\} \\
 & + \frac {2q} {\sqrt{h}} \sin \bigg( \frac{\alpha}{2}\bigg) \, \log\bigg( \frac{h+1+2\sqrt{h}\,\cos(\alpha/2)}{\sqrt{1+h^2-2h\cos \alpha}} \bigg).
\end{align*}
Inserting all of the above auxiliary calculations into expression $(\ref{mhsf})$, after simplifications, we obtain the following expression for $\cF(C_{S,\alpha})$,
\begin{equation*}
\cF(C_{S,\alpha}) = \frac {\pi} {\pi-\alpha+\sin \alpha} \bigg\{ 1 + \frac {q(h+1)}{2h} \bigg(1 - \frac {\alpha} {\pi} \bigg) - \frac {q(h-1)}{\pi h} \tan^{-1}\bigg( \cot\bigg(\frac{\alpha}{2}\bigg) \, \frac {h-1} {h+1} \bigg) \bigg\}.
\end{equation*}
Hence the functional $\cF(C_{S,\alpha})$ assumes the claimed form.

\qed


\noindent{\bf Proof of Proposition \ref{polyextf}.} Inserting expression $(\ref{polyextf})$ into the formula for the Mhaskar-Saff functional $(\ref{mhsf})$, we see that
\begin{empheq}{align*} 
\cF(C_{S,\alpha})  = \frac {\pi} {2(\pi-\alpha+\sin \alpha)} \bigg\{2  & + \int_\alpha^\pi (a\cos^2\varphi+b\sin\varphi+c) \sin\varphi\, d\varphi \\ 
							& + \frac {2} {\pi} \int_\alpha^\pi (a\cos^2\varphi+b\sin\varphi+c)  \sqrt{ \frac {1-\cos \alpha} {\cos \alpha - \cos\varphi} } \sin\varphi\, d\varphi \\ 
			& - \frac {2} {\pi} \int_\alpha^\pi (a\cos^2\varphi+b\sin\varphi+c)  \tan^{-1} \sqrt{ \frac {1-\cos \alpha} {\cos \alpha - \cos\varphi} } \sin\varphi\, d\varphi \bigg\}.
\end{empheq}
Each of the integrals on the right hand side of the last expression, can be easily evaluated via the substitution $u=\cos\varphi$. Evaluating those integrals, and inserting them into
the expression for $\cF(C_{S,\alpha})$, after some algebra, we arrive at the desired expression $(\ref{msfexpquadpoly})$.

\qed

\noindent{\bf Proof of Theorem \ref{transceqsupppolyextf}.} Let 
\begin{empheq}{align}\label{auxmsfq1}
 w(\alpha)   =  \frac {1} {36(\pi-\alpha+\sin \alpha)} &  \bigg\{  \tan \left( \frac{\alpha}{2}\right) \bigg [32\, a \cos^3 \alpha  + 4(2a+9b) \cos^2 \alpha \\  \nonumber
  + 4(9c-5a) \cos \alpha & + 4 a - 36 b + 36 c \bigg]  + 12 (a+3c) (\pi-\alpha) + 36 \pi \bigg\}.
\end{empheq}
Differentiating $(\ref{auxmsfq1})$ with respect to $\alpha$ and simplifying, we obtain
\begin{equation}\label{dfw}
w'(\alpha) = \frac {1} {3(\pi - \alpha + \sin \alpha)} \, \omega(\alpha),
\end{equation}
where
\begin{equation}\label{omegaqef}
\omega(\alpha) = 6 \sin^2\bigg(\frac{\alpha}{2} \bigg) w(\alpha) + 8a\cos^3\alpha + 2(3b-2a)\cos^2\alpha + (3c-5a-3b)\cos\alpha + a - 3b -3c. 
\end{equation}
Therefore, the critical points of $w(\alpha)$ are given by solutions of the transcendental equation
\[
6 \sin^2\bigg(\frac{\alpha}{2} \bigg) w(\alpha) + 8a\cos^3\alpha + 2(3b-2a)\cos^2\alpha + (3c-5a-3b)\cos\alpha + a - 3b -3c = 0,
\]
or, equivalently,
\begin{empheq}{align*}
 8a\cos^3 \alpha (2\sin\alpha+3(\pi-\alpha)) & + \cos^2\alpha((2a+9b)\sin\alpha - 6(2a-3b)(\pi-\alpha)) \\  
& + \frac{1}{2}\sin(2\alpha) (9b-22a) + 3(2a-3b)(\pi-\alpha)+9\pi \\   
& - 9\cos\alpha (\pi + (2a+b)(\pi-\alpha)) - 2\sin\alpha(2a-9b) =0.
\end{empheq}
Rearranging the latter, we obtain $(\ref{alsupppoly})$.

We proceed by showing the existence and uniqueness of a critical point of $w(\alpha)$.  We will be making use of an argument developed in \cite{bds1}. First, observe that
\begin{align*}
\lim_{\alpha \rightarrow \pi-} \omega (\alpha)  & = 6\lim_{\alpha \rightarrow \pi-} w(\alpha) - 6a + 6b -6c \\
								 & = 6 \pi \lim_{\alpha \rightarrow \pi-} \frac {1} {\pi - \alpha + \sin \alpha} - 6a + 6b -6c \\
								& =  +\infty.
\end{align*}
Hence, there is a smallest $\alpha_0\in[0,\pi)$ such that $\omega(\alpha)>0$ for $\alpha\in(\alpha_0,\pi)$. If $\alpha_0=0$, then $w(\alpha)$ is 
strictly increasing on $(0,\pi)$, and attains minimum at $\alpha=0$. If $\alpha_0>0$, we have that $\omega(\alpha)>0$ for $\alpha\in(\alpha_0,\pi)$. Taking into account the continuity of $\omega(\alpha)$, by passing to the limit $\alpha\rightarrow \alpha_0+$ in the latter inequality, we infer that $\omega(\alpha_0)\geq 0$. Since $\alpha_0$ was the smallest $\alpha$ such that $\omega(\alpha)>0$ on $(\alpha_0,\pi)$, we deduce that $\omega(\alpha_0) = 0$.

As the sign of $w'(\alpha)$ is determined by the sign of $\omega(\alpha)$, it is clear that $w'(\alpha)>0$ on $(\alpha_0,\pi)$, and $w'(\alpha_0)=0$. Next, suppose that $\xi\in(0,\pi)$ is a critical point of $w(\alpha)$, that is $w'(\xi)=0$. Using expression $(\ref{dfw})$, we easily find that
\begin{empheq}{align}\label{2ndder}
w''(\alpha) & = \frac {1} {3(\pi-\alpha+\sin \alpha)} \bigg\{  12\sin^2\bigg(\frac{\alpha}{2}\bigg) w'(\alpha)  \\ \nonumber 
	& + \sin(\alpha)\big(3w(\alpha)-24a\cos^2\alpha-4(3b-2a)\cos\alpha -(3c-5a-3b)\big) \bigg\}.
\end{empheq}
Our goal is to show that $w''(\xi)>0$. Indeed, from $(\ref{2ndder})$ it is clear that
\[
w''(\xi) =  \frac {\sin(\xi)} {3(\pi-\xi+\sin \xi)} \bigg( 3w(\xi)-24a\cos^2\xi-4(3b-2a)\cos\xi -(3c-5a-3b) \bigg),
\]
and so the sign of $w''(\xi)$ is determined by the sign of the expression $3w(\xi)-24a\cos^2\xi-4(3b-2a)\cos\xi -(3c-5a-3b)$. As $\xi$ is a critical point of $w$, from $(\ref{dfw})$ and $(\ref{omegaqef})$ it follows that 
\[
3(1-\cos\xi)w(\xi) = -\big(8a\cos^3\xi+2(3b-2a)\cos^2\xi+(3c-5a-3b)\cos\xi+a-3b-3c \big).
\]
Taking into account that the sign of $w''(\xi)$ is the same as the sign of $w''(\xi)(1-\cos\xi)$, it is a straightforward calculation to see that the sign of $w''(\xi)$ is, in fact, determined by the sign of 
\begin{equation*}
m(\xi) = 16\cos^3\xi+2(3\widetilde{b}-14)\cos^2\xi+4(2-3\widetilde{b})\cos\xi+4+6\widetilde{b},
\end{equation*}
where we denoted $\widetilde{b}:=b/a>2$.
Setting $t=\cos\xi\in(-1,1)$, abusing notation, we consider the function
\begin{equation*}
m(t) = 16t^3+2(3\widetilde{b}-14)t^2+4(2-3\widetilde{b})t+4+6\widetilde{b}, \quad -1< t < 1.
\end{equation*}
First assume that $0\leq t <1$. Then $m'(t)<0$ on $[0,1)$, so that $(t)$ is strictly decreasing function on $[0,1)$. This obviously entails $m(t)>m(1)$, for all $t\in[0,1)$. A direct calculation shows that $m(1)=0$, and thus $m(t)>0$ for $0\leq t <1$.
If $-1< t \leq 0$, it follows that we need to demonstrate the positivity of 
\[
p(t) = -16t^3+2(3\widetilde{b}-14)t^2-4(2-3\widetilde{b})t+4+6\widetilde{b}, \quad 0 \leq t < 1.
\]
Using the trivial estimates, we see that $p(t)>0$, for $0 \leq t < 1$. Combining the previous two estimates, we deduce that $m(\xi)>0$. Thus $w''(\xi)  > 0$, as desired.

We are now able to conclude that any critical point of $w(\alpha)$ is a minimum of $w(\alpha)$. Since $w(\alpha)$ is smooth on $(0,\pi)$, we see that $\alpha_0$ is a unique local minimum for $w(\alpha)$ on $(0,\pi)$. Thus we are only left to show that $\alpha_0$ is a global minimum as well. Indeed, we observe that $w'(\alpha)>0$ on $(\alpha_0,\pi)$, and $w'(\alpha)<0$ on $(0,\alpha_0)$. Hence $\omega(\alpha)>0$ on $(\alpha_0,\pi)$, and $\omega(\alpha)<0$ on $(0,\alpha_0)$. This means that $w(\alpha)$ has exactly one global minimum on $[0,\pi)$, which is either a unique solution $\alpha_0\in(0,\pi)$ of the transcendental equation $(\ref{alsupppoly})$, if it exists, or $\alpha_0=0$, if such a solution fails to exist.

\qed


\noindent{\bf Proof of Theorem \ref{transceqsupppolyextf}.} Inserting $Q(\varphi)=a\cos^2\varphi+b\cos\varphi+c$ into $(\ref{auxss})$, we easily find
\begin{equation}\label{gtpoly}
g(t) = -\frac{\sin t}{12\pi \sqrt{1+\cos t}} (8a\cos^2t + 2(2a+3b)\cos t - a + 3b + 3c ).
\end{equation}
Next, according to $(\ref{auxfs})$, we need to evaluate 
\begin{equation*}
I(\varphi) = \int_{\alpha_0}^{\varphi} \frac{ g(t) \sin t \,dt} {\sqrt{\cos t - \cos\varphi}}.
\end{equation*}
Substituting $(\ref{gtpoly})$ into the above integral, we see that
\begin{align*}
I(\varphi) = - \frac {1} {12\pi} \bigg\{ 8a  & \int_{\alpha_0}^{\varphi} \frac{ \cos^2 t \sin^2 t \,dt} {\sqrt{1+\cos t} \sqrt{\cos t - \cos\varphi}} \\
			+  2(2a+3b)  & \int_{\alpha_0}^{\varphi} \frac{ \cos t \sin^2 t \,dt} {\sqrt{1+\cos t} \sqrt{\cos t - \cos\varphi}}  \\
			+ (3b + 3c - a) & \int_{\alpha_0}^{\varphi} \frac{ \sin^2 t \,dt} {\sqrt{1+\cos t} \sqrt{\cos t - \cos\varphi}} \bigg\}.
\end{align*}
All of the integrals on the right hand side of the last expression are readily evaluated via the elementary substitution $u=\cos t$. Evaluating the integrals and simplifying, we see that
\begin{empheq}{align}\label{uaxicalc}
I(\varphi)  = & - \frac {1} {12\pi} \bigg[ 8a \bigg\{ \frac{1}{48} \sqrt{1-\cos\alpha_0} \sqrt{\cos\alpha_0-\cos\varphi} (16\cos^2\alpha_0 \\ \nonumber
			& + 4(5\cos\alpha_0-2) \cos\varphi) -4\cos\alpha_0 + 15\cos(2\varphi)+9) \\ \nonumber
			& + \frac{1}{16} (1-\cos\varphi) (4\cos\varphi + 5\cos(2\varphi)+7)\tan^{-1}\sqrt{\frac{\cos\alpha_0 - \cos\varphi}{1-\cos\alpha_0}} \bigg\} \\ \nonumber
			& + 2(2a+3b) \bigg\{ \frac{1}{4}  \sqrt{1-\cos\alpha_0} \sqrt{\cos\alpha_0-\cos\varphi} (2\cos\alpha_0 + 3\cos\varphi-1) \\ \nonumber
			& + \frac{1}{4}  (1-\cos\varphi) (3\cos\varphi+1) \tan^{-1}\sqrt{\frac{\cos\alpha_0 - \cos\varphi}{1-\cos\alpha_0}} \bigg\} \\ \nonumber
			& + (3b+3c-a) \bigg\{  (1-\cos\varphi) \tan^{-1}\sqrt{\frac{\cos\alpha_0 - \cos\varphi}{1-\cos\alpha_0}} \\ \nonumber
			& +  \sqrt{1-\cos\alpha_0} \sqrt{\cos\alpha_0-\cos\varphi} \bigg\} \bigg].
\end{empheq}
Differentiating the last expression, and upon inserting the result into $(\ref{auxfs})$, after some algebra, we deduce that
\begin{empheq}{align*}
F(\varphi)  = & \frac {1}{36\,\pi^2} \bigg\{  \sqrt{1-\cos\alpha_0}\, \sqrt{\cos\alpha_0-\cos\varphi}\, (20a\cos\alpha_0 + 60a\cos\varphi + 10a + 27b) \\
 & - \sqrt{\frac{1-\cos\alpha_0}{\cos\alpha_0-\cos\varphi}}\, (8a\cos^2\alpha_0 + 10a\cos\alpha_0 \cos\varphi + (4a+9b) \cos\alpha_0 + (20a + 27b) \cos\varphi \\ 
		& + 15a\cos(2\varphi) + 9a + 18b + 18c) \\
		& - 6 \tan^{-1} \sqrt{\frac{\cos\alpha_0-\cos\varphi}{1-\cos\alpha_0}}\, (15a\cos^2\varphi + 9b\cos\varphi - 4a + 3c) \bigg\}.
\end{empheq}
At last, we compute the Robin constant $F_Q$. By $(\ref{dcfs})$ and $(\ref{auxfs})$, taking into account that
\begin{align*}
F(\varphi) & = \frac {2} {\pi} \frac {1} {\sin\varphi} \, \frac {d} {d\varphi} \int_{\alpha_0}^{\varphi} \frac{ g(t) \sin t \,dt} {\sqrt{\cos t - \cos\varphi}} \\
		& = \frac {2} {\pi} I'(\varphi),
\end{align*}
we easily see that formula $(\ref{dcfs})$ reduces to
\[
F_Q = \frac {\pi} {\sin\alpha_0 + \pi - \alpha_0} \bigg\{ 1 + 4 (I(\alpha_0) - I(\pi)) \bigg\}.
\]
Inserting $(\ref{uaxicalc})$ into the latter expression, after simplifications, we infer that
\begin{empheq}{align*}
F_Q & =  \frac {1} {36(\pi-\alpha_0+\sin \alpha_0)}  \bigg\{  \tan \left( \frac{\alpha_0}{2}\right) \bigg [32\, a \cos^3 \alpha_0  + 4(2a+9b) \cos^2 \alpha_0 \\  \nonumber
& + 4(9c-5a) \cos \alpha_0 + 4 a - 36 b + 36 c \bigg]  + 12 (a+3c) (\pi-\alpha_0) + 36 \pi \bigg\}.
\end{empheq}
This completes the proof of Theorem \ref{transceqsupppolyextf}.
\qed

\section{Acknowledgements} The author would like to thank his doctoral advisor Prof. Igor E. Pritsker for suggesting 
the problem and stimulating discussions. The author expresses gratitude to Prof. Edward B. Saff for his kind permission to
use a reference from the forthcoming book \cite{bhs}.

\end{document}